\documentclass{amsart}
\usepackage{amsmath, amssymb, latexsym, amsthm, mathrsfs}
\usepackage[all]{xy}

\newcommand{\ed}{

\subsection*{Acknowledgments} This work is an extension of a part of the first named author's M.Sc.\ thesis
at the Weizmann Institute of Science, supervised by Gady Kozma and the second named author.
We thank Gady Kozma for useful discussions, and the Weizmann Institute of Science for the stimulating atmosphere.
We owe special thanks to Lyubomyr Zdomskyy for reading the paper and making useful comments.

\end{document}}

      \newenvironment{changemargin}[2]{\begin{list}{}{
         \setlength{\topsep}{0pt}\setlength{\leftmargin}{0pt}
         \setlength{\rightmargin}{0pt}
         \setlength{\listparindent}{\parindent}
         \setlength{\itemindent}{\parindent}
         \setlength{\parsep}{0pt plus 1pt}
         \addtolength{\leftmargin}{#1}\addtolength{\rightmargin}{#2}
         }\item }{\end{list}}

\newcommand{\Pa}[8]{\bibitem{#1} {#2}, \emph{#3}, {#4} \textbf{#5} ({#6}), {#7}--{#8}.}

\newcommand{\Bc}[9]{\bibitem{#1} {#2}, \emph{#3}, in: \textbf{#4} (#5), #6 #7, #8--#9.}

\newcommand{\Setting}[7]{\xymatrix@R=4pt@C=7pt{#1\ar@{-}[r]&#2\ar@{-}[r]&#3\\&#4\ar@{-}[u]\\
#5\ar@{-}[uu]\ar@{-}[r] & #6\ar@{-}[u]\ar@{-}[r] & #7\ar@{-}[uu]}}

\newcommand{\Bgp}{{\Z^\N}}

\newcommand{\arx}[1]{\texttt{http://arxiv.org/math/#1}}
\newcommand{\bq}{\begin{quote}}
\newcommand{\eq}{\end{quote}}

\newcommand{\CH}{the Continuum Hypothesis}

\newcommand{\inv}{^{-1}}
\newcommand{\Cantor}{{\{0,1\}^\N}}

\newcommand{\sr}[2]{{\txt{$#1$\\#2}}}
\newcommand{\nsr}[2]{#1}
\newcommand{\N}{\mathbb{N}}
\newcommand{\NN}{{\N^{\N}}}

\newcommand{\roth}{{[\N]^{\oo}}}%{{[\N]^{\aleph_0}}}
\newcommand{\Fin}{{[\N]^{<\oo}}}%{{[\N]^{<\aleph_0}}}
\newcommand{\seq}[1]{\{#1\}_{n\in\N}}

\newcommand{\sseq}[1]{\{#1 : n\in\N\}}

\newcommand{\op}{\operatorname}

\newcommand{\scrA}{\mathscr{A}}
\newcommand{\scrB}{\mathscr{B}}

\newcommand{\B}{\mathrm{B}}

\newcommand{\BG}{\B_\Gamma}

\newcommand{\BO}{\B_\Omega}

\newcommand{\cF}{\mathcal{F}}

\newcommand{\cM}{\mathcal{M}}

\newcommand{\rmO}{\mathrm{O}}

\newcommand{\R}{\mathbb{R}}
\newcommand{\cU}{\mathcal{U}}
\newcommand{\Union}{\bigcup}
\newcommand{\cV}{\mathcal{V}}
\newcommand{\cW}{\mathcal{W}}
\newcommand{\Z}{{\mathbb Z}}
\newcommand{\Impl}{\Rightarrow}
\long\def\forget#1\forgotten{}
\newcommand{\ft}{\mathfrak{t}}
\newcommand{\fb}{\mathfrak{b}}
\newcommand{\fc}{\mathfrak{c}}
\newcommand{\fd}{\mathfrak{d}}

\newcommand{\oo}{\infty}

\newcommand{\fp}{\mathfrak{p}}

\newcommand{\x}{\times}

\newcommand{\nin}{\notin}

\newcommand{\sbst}{\subseteq}

\newcommand{\sm}{\setminus}
\newcommand{\as}{\subseteq^*}%{\let\proclaim\relax}

\renewcommand{\>}{\rangle}

\newcommand{\cov}{\op{cov}}

\newcommand{\non}{\op{non}}

\newtheorem{thm}{Theorem}[section]
\newcommand{\bthm}{\begin{thm}} \newcommand{\ethm}{\end{thm}}
\newtheorem{prop}[thm]{Proposition}
\newcommand{\bprp}{\begin{prop}} \newcommand{\eprp}{\end{prop}}
\newtheorem{fact}[thm]{Fact}
\newcommand{\bfct}{\begin{fact}} \newcommand{\efct}{\end{fact}}
\newtheorem{prob}[thm]{Problem}
\newcommand{\bprb}{\begin{prob}} \newcommand{\eprb}{\end{prob}}
\newtheorem{lem}[thm]{Lemma}
\newcommand{\blem}{\begin{lem}} \newcommand{\elem}{\end{lem}}
\newtheorem{claim}[thm]{Claim}
\newcommand{\bclm}{\begin{claim}} \newcommand{\eclm}{\end{claim}}
\newtheorem{cor}[thm]{Corollary}
\newcommand{\bcor}{\begin{cor}} \newcommand{\ecor}{\end{cor}}
\newtheorem{conj}[thm]{Conjecture}
\newcommand{\bcnj}{\begin{conj}} \newcommand{\ecnj}{\end{conj}}
\theoremstyle{definition}
\newtheorem{defn}[thm]{Definition}
\newcommand{\bdfn}{\begin{defn}} \newcommand{\edfn}{\end{defn}}
\theoremstyle{remark}
\newtheorem{rem}[thm]{Remark}
\newcommand{\brem}{\begin{rem}} \newcommand{\erem}{\end{rem}}
\newtheorem{cnv}[thm]{Convention}
\newcommand{\bcnv}{\begin{cnv}} \newcommand{\ecnv}{\end{cnv}}
\newtheorem{exam}[thm]{Example}
\newcommand{\bexm}{\begin{exam}} \newcommand{\eexm}{\end{exam}}
\newcommand{\bpf}{\begin{proof}} \newcommand{\epf}{\end{proof}}
\newcommand{\be}{\begin{enumerate}}
\newcommand{\ee}{\end{enumerate}}
\newcommand{\bi}{\begin{itemize}}
\newcommand{\itm}{\item}
\newcommand{\ei}{\end{itemize}}
\forget
\forgotten

\newcommand{\sone}{\mathsf{S}_1}
\newcommand{\sfin}{\mathsf{S}_\mathrm{fin}}
\newcommand{\ufin}{\mathsf{U}_\mathrm{fin}}

\title[Linear $\sigma$-additivity]{Linear $\sigma$-additivity and some applications}
\author{Tal Orenshtein}
\address[Orenshtein]{Department of Mathematics,
Weizmann Institute of Science, Rehovot 76100, Israel}
\email{talo@weizmann.ac.il}
\urladdr{http://www.orenshtein.com}
\author{Boaz Tsaban}
\address[Tsaban]{Department of Mathematics, Bar-Ilan University, Ramat-Gan 52900, Israel}
\email{tsaban@math.biu.ac.il}
\urladdr{http://www.cs.biu.ac.il/\~{}tsaban}

\begin{document}

\begin{abstract}
We show that countable
increasing unions preserve a large family of well-studied covering
properties, which are not necessarily $\sigma$-additive.
Using this, together with infinite-combinatorial methods and simple forcing theoretic
methods, we explain several phenomena, settle problems of
Just, Miller, Scheepers and Szeptycki \cite{coc2},
Gruenhage and Szeptycki \cite{FUfin},
Tsaban and Zdomskyy \cite{SFT},
and Tsaban \cite{o-bdd, OPiT},
and construct topological groups with very strong combinatorial properties.
\end{abstract}

\maketitle

\section{Introduction}

The following natural definition unifies all results presented here.

\bdfn\label{LSA}
Let $\cF$ be a family of topological spaces. $\cF$ is \emph{linearly $\sigma$-additive}
if it is preserved by countable increasing unions. That is:
For each topological space $X=\Union_n X_n$ with $X_1\sbst X_2\sbst\dots$ and $X_n\in\cF$ for
all $n\in\N$, $X\in\cF$.
\edfn

Removing the restriction that $X_n\sbst X_{n+1}$ for all $n$, we obtain
the definition of \emph{$\sigma$-additive} family.
We identify a topological property with the family of all topological spaces satisfying it.
Thus, we may talk about linearly $\sigma$-additive \emph{properties}.

We consider additivity in the context of topological selection
principles, to which we now give a brief introduction.\footnote{Extended introductions
to this field are available in \cite{KocSurv, LecceSurvey, ict}.}
This is a framework suggested by Scheepers in \cite{coc1}
to study in a uniform manner a variety of properties introduced
in different mathematical disciplines, since the early
1920's, by Menger, Hurewicz, Rothberger, and Gerlits and Nagy, and many others.

Let $X$ be a topological space. We say that $\cU$ is a \emph{cover}
of $X$ if $X=\Union\cU$, but $X\nin\cU$.
Often, $X$ is considered as a subspace of another space $Y$,
and in this case we always consider covers of $X$ by subsets of $Y$,
and require instead that no member of the cover contains $X$.
Let $\rmO(X)$ be the family of
all open covers of $X$. Define the following subfamilies of $\rmO(X)$:
$\cU\in\Omega(X)$ if each finite subset of $X$ is contained in some member of $\cU$.
$\cU\in\Gamma(X)$ if $\cU$ is infinite, and each element of $X$ is contained in all but
finitely many members of $\cU$.

Some of the following statements may hold for families $\scrA$ and $\scrB$ of covers of $X$.
\begin{description}
\item[$\sone(\scrA,\scrB)$] For all $\cU_1,\cU_2,\dots\in\scrA$, there are
$U_1\in\cU_1,U_2\in\cU_2,\dots$ such that $\sseq{U_n}\in\scrB$.
\item[$\sfin(\scrA,\scrB)$] For all $\cU_1,\cU_2,\dots\in\scrA$, there are
finite $\cF_1\sbst\cU_1,\cF_2\sbst\cU_2,\dots$ such that $\Union_n\cF_n\in\scrB$.
\item[$\ufin(\scrA,\scrB)$] For all $\cU_1,\cU_2,\dots\in\scrA$, none containing
a finite subcover, there are finite $\cF_1\sbst\cU_1,\cF_2\sbst\cU_2,\dots$ such that $\sseq{\Union\cF_n}\in\scrB$.
\end{description}
We say, e.g., that $X$ satisfies $\sone(\rmO,\rmO)$ if the statement $\sone(\rmO(X),\rmO(X))$ holds.
This way, $\sone(\rmO,\rmO)$ is a property of topological spaces, and similarly for all other statements
and families of covers.
Under some mild hypotheses on the considered topological spaces,
each nontrivial property among these properties, where $\scrA,\scrB$ range over $\rmO,\Omega,\Gamma$,
is equivalent to one in Figure \ref{SchDiag}  \cite{coc1, coc2}. In this diagram, an arrow denotes implication.

\begin{figure}[!htp]
\begin{changemargin}{-4cm}{-3cm}
\begin{center}
{%\scriptsize
$\xymatrix@R=8pt{
%1
&
&
& \sr{\ufin(\rmO,\Gamma)}{Hurewicz}\ar[r]
& \sr{\ufin(\rmO,\Omega)}{}\ar[rr]
& & \sr{\sfin(\rmO,\rmO)}{Menger}
\\
%2
&
&
& \sr{\sfin(\Gamma,\Omega)}{}\ar[ur]
\\
%3
& \sr{\sone(\Gamma,\Gamma)}{}\ar[r]\ar[uurr]
& \sr{\sone(\Gamma,\Omega)}{}\ar[rr]\ar[ur]
& & \sr{\sone(\Gamma,\rmO)}{}\ar[uurr]
\\
%4
&
&
& \sr{\sfin(\Omega,\Omega)}{Menger$^\uparrow$}\ar'[u][uu]
\\
%5
& \sr{\sone(\Omega,\Gamma)}{Gerlits-Nagy}\ar[r]\ar[uu]
& \sr{\sone(\Omega,\Omega)}{Rothberger$^\uparrow$}\ar[uu]\ar[rr]\ar[ur]
& & \sr{\sone(\rmO,\rmO)}{Rothberger}\ar[uu]
}$
}
\caption{The Scheepers Diagram}\label{SchDiag}
\end{center}
\end{changemargin}
\end{figure}

In this diagram, the classical name of a property is indicated below it,
as well as two names ending with a symbol $\uparrow$, by which
we indicate that the properties $\sone(\Omega,\Omega)$ and
$\sfin(\Omega,\Omega)$ may also be viewed as classical ones \cite{Sakai88, coc2},
in accordance with the following notation.

\bdfn\label{pows}
Let $P$ be a property of topological spaces. \emph{$X$ satisfies $P^\uparrow$} if
all finite powers $X^k$ of $X$ satisfy $P$.
\edfn

The Scheepers diagram is at the heart of the field of topological selection principles,
and many additional---classical and new---properties are studied in relation to it.
The reader is encouraged to consider his favorite properties in light of the
results presented here.

\medskip

In Section \ref{AllLSAsec}, we prove that all properties in the Scheepers Diagram
are linearly $\sigma$-additive, and are thus hereditary for $F_\sigma$ subsets.
This solves a problem of Tsaban and Zdomskyy from \cite{SFT}.

A crucial part of the proof that the studied properties are linearly $\sigma$-additive,
is a recent theorem of F. Jordan \cite{Jordan}.
Miller asked in \cite{MilRelGa} whether $\sone(\Omega,\Gamma)$ is linearly $\sigma$-additive.
A negative solution would have solved the notorious Gerlits-Nagy problem \cite{GN}.
Using a brilliant argument, Jordan proved that this is not the case.
We give a direct version of Jordan's solution, and provide some applications.

In Section \ref{JMSSSec}, we use Jordan's method in a proof that if there is an unbounded family
of cardinality $\aleph_1$ in the Baire space, then there is an uncountable set
of real numbers, satisfying $\sone(\Omega,\Gamma)$. This settles in the positive a problem
from the seminal paper of Just, Miller, Scheepers, and Szeptycki \cite{coc2}.
Indeed, our result is more general, and also solves a problem of Gruenhage and Szeptycki \cite{FUfin}.

In Section \ref{heredsec}, we apply linear $\sigma$-additivity to study heredity of properties,
answer a question of Zdomskyy, and suggest a simple revision
of a question of Bukovsk\'y, Rec\l{}aw, and Repick\'y, which makes it possible to answer
it in the positive. (The problem, as originally stated, was answered in the negative by Miller \cite{MilNonGamma}.)
This section also explains the phenomenon observed in \cite{ideals}, that
none of the considered properties is hereditary in the open case, whereas in the Borel case
some are and some are not hereditary.

In Section \ref{topgpssec} we apply our results to construct
topological groups with strong combinatorial properties, and solve
problems from \cite{o-bdd} and \cite{OPiT}.

\medskip

\subsection{Generalizations}
The results presented here require little or no assumptions on the topology
of the studied spaces. However, they are interesting even when restricting
attention to, e.g., metric spaces or even subsets of $\R$.

For concreteness, we present the results only for the types of covers mentioned
above, but the proofs show that they hold for many additional types.
In particular, define $\B$, $\BO$, $\BG$ as $\rmO$, $\Omega$, $\Gamma$ were
defined, replacing \emph{open cover} by \emph{countable Borel cover}.
The properties thus obtained have rich history of their own \cite{CBC}.
Mild assumptions on $X$ imply that the considered open covers may be
assumed to be countable, and this makes the Borel variants of the studied properties
(strictly) stronger \cite{CBC}.
All of the results presented here also hold in the Borel case
(after replacing \emph{open} or \emph{closed} by \emph{Borel}).
Moreover, unlike some of the results in the open case,
none of the Borel variant requires any assumption on the topology of $X$.

\section{Linear $\sigma$-additivity in the Scheepers Diagram}\label{AllLSAsec}

One motivation for studying linear $\sigma$-additivity in the context
of the Scheepers Diagram is an experimentally observed dichotomy concerning additivity in of properties
in this diagram:
Each property there is either provably $\sigma$-additive, or else not even provably \emph{finitely} additive.\footnote{A
survey of the involved results, with complete proofs, is available in \cite{AddQuad}.}
In this section, we prove the following.

\bthm\label{AllLSA}
All properties in the Scheepers Diagram (Figure \ref{SchDiag}) are linearly $\sigma$-additive.
\ethm

It remains to notice the following.

\blem\label{dicho}
Each linearly $\sigma$-additive property is either $\sigma$-additive, or else
not finitely additive.
\elem
\bpf
$\Union_nX_n=\Union_n(\Union_{m\le n}X_m)$.
\epf

\bcor
Each property in the Scheepers Diagram is either $\sigma$-additive, or else
not additive.\qed
\ecor

Before proving Theorem \ref{AllLSA}, we point out two additional consequences.
In Problem 4.9 of \cite{SFT} and Problem 6.2 of \cite{OPiT}, Tsaban and Zdomskyy ask whether
$\sone(\Gamma,\Omega)$ and $\sfin(\Gamma,\Omega)$ are hereditary for $F_\sigma$ subsets.
We obtain a positive answer.

\bcor\label{Fs}
All properties in the Scheepers Diagram are hereditary for $F_\sigma$ subsets.
\ecor
\bpf
These properties are easily seen to be hereditary for closed subsets \cite{coc2},
and countable unions of closed subsets can be presented as countable increasing
unions of closed subsets. Apply Theorem \ref{AllLSA}.
\epf

The following was observed in the past for at least some properties in the Scheepers
Diagram, each time using an ad-hoc argument.

\bcor
Let $P$ be a property in the Scheepers Diagram, and $X,D$ be subspaces of some topological space.
If $X$ satisfies $P$ and $D$ is countable, then $X\cup D$ satisfies $P$.
\ecor
\bpf
In light of Theorem \ref{AllLSA}, it remains to observe that for each singleton $\{a\}$,
$X\cup\{a\}$ satisfies $P$. This is not hard to verify.
\epf

Theorem \ref{AllLSA} is proved in parts. The properties
$\sone(\rmO,\rmO)$, $\sone(\Gamma,\rmO)$, $\sfin(\rmO,\rmO)$, $\sone(\Gamma,\Gamma)$, and $\ufin(\rmO,\Gamma)$
are in fact $\sigma$-additive \cite{coc2, wqn} (see \cite{AddQuad}).

Linear $\sigma$-additivity of $\ufin(\rmO,\Omega)$ was proved in \cite{SFT} for sets of reals.
It also follows from the following.

\bthm\label{AOm}
For all $\Pi\in\{\sone,\sfin,\ufin\}$ and $\scrA\in\{\Gamma,\Omega,\rmO\}$,
$\Pi(\scrA,\Omega)$ is linearly $\sigma$-additive.
\ethm
\bpf
We prove the theorem for $\Pi=\sfin$; the other proofs being similar.

Assume that $X=\Union_nX_n$ is an increasing union, with each
$X_n$ satisfying $\sfin(\scrA,\Omega)$.
Let $\cU_1,\cU_2,\dots\in\scrA(X)$.
We first exclude the trivial case: Assume that for infinitely many
$n$, there are $m_n$ and elements $U_{m_n}\in\cU_{m_n}$ such that $X_n\sbst U_{m_n}$.
As $X$ is not contained in any member of any $\cU_n$ and the sets $X_n$ increase to $X$,
we may if necessary thin out the sequence $m_n$ to make it increasing.
Then $\sseq{U_{m_n}}\in\Omega(X)$, and this suffices.

Thus, we may assume that for all $n,k$, we have that $\cU_k\in\scrA(X_n)$.
Take a partition $\N=\Union_nI_n$ of $\N$ into infinite sets $I_n$.
Fix $n$.
As $X_n$ satisfies $\sfin(\scrA,\Omega)$, there are finite $\cF_k\sbst\cU_k$,
$k\in I_n$, such that $\Union_{k\in I_n}\cF_k\in\Omega(X_n)$.
Then $\Union_{k\in\N}\cF_k\in\Omega(X)$.
\epf

The remaining property, $\sone(\Omega,\Gamma)$, was treated by F. Jordan.

\subsection{Jordan's Theorem and some applications}

%An infinite $\cF\sbst P(X)$ is \emph{point-cofinite} for $Y\sbst X$ if
%for each $y\in Y$, $\{A\in\cF : y\in A\}$ is a cofinite subset of $\cF$.
The following technical lemma will be useful in the proof of Jordan's Theorem below.

\blem\label{capinfgg}
Let $Y\sbst X$ be such that $Y$ satisfies $\sone(\Gamma,\Gamma)$.
Assume that for each $n$,
\be
\itm $\cU_n$ is an infinite family of open subsets of $X$; and
\itm For each $y\in Y$, $y\in U$ for all but finitely many $U\in\cU_n$.
\ee
Then there are infinite $\cV_1\sbst\cU_1, \cV_2\sbst\cU_2,\dots$, such that
for each $y\in Y$, $y\in\bigcap\cV_n$ for all but finitely many $n$.
\elem
\bpf
It may be the case than no subset of $\cU_n$ is in $\Gamma(Y)$.

Case 1: For all but finitely many $n$, $\cV_n=\{U\in\cU_n : Y\sbst U\}$
is infinite. Then the sets $\cV_n$ thus defined are as required.

Case 2: Let $I$ be the set of all $n$ such that $\cV_n=\{U\in\cU_n : Y\sbst U\}$
is infinite, and $J=\N\sm I$. For each $n\in J$, $\cW_n=\{U\in\cU_n : Y\not\sbst U\}$ is infinite,
and thus $\cW_n\in\Gamma(Y)$. As $J$ is infinite, there are by Theorem 15 of \cite{capinf} infinite
$\cV_n\sbst\cW_n$, $n\in J$, such that each $y\in Y$ belongs to $\bigcap\cV_n$ for all but finitely
many $n$. (Briefly: By thinning out if needed, we may assume that each $\cW_n$ is countable
and that $\cW_n\cap\cW_m=\emptyset$ for
$m\neq n$ \cite{coc1}. Apply $\sone(\Gamma,\Gamma)$ to the countable family of all cofinite subsets of all $\cW_n$
to obtain $\cV\in\Gamma(Y)$.
Let $\cV_n=\cV\cap\cW_n$.)
The sets $\cV_1,\cV_2,\dots$ are as required.
\epf

\bthm[Jordan \cite{Jordan}]
$\sone(\Omega,\Gamma)$ is linearly $\sigma$-additive.
\ethm
\bpf
We give a direct proof, following what seems to be the essence of Jordan's arguments.
The following statement can be deduced from Theorem 7 in \cite{Jordan}.
\blem\label{jlem}
Let $X=\Union_nX_n$ be an increasing union, where each $X_n$ satisfies
$\sone(\Gamma,\Gamma)$. For all $\cU_1\in\Gamma(X_1),\cU_2\in\Gamma(X_2),\dots$, there are
infinite $\cV_1\sbst\cU_1, \cV_2\sbst\cU_2,\dots$, such that for each $x\in X$, $x\in \bigcap\cV_n$
for all but finitely many $n$.
\elem
\bpf
Step 1: By Lemma \ref{capinfgg}, we may thin out the families $\cU_n$, so that
they remain infinite, and each member of $X_1$ belongs to all but finitely many $\bigcap\cU_n$.
Let $\cV_1=\cU_1$.

Step 2: By the same Lemma, we may thin out further the families $\cU_n$, $n\ge 2$, so that
they remain infinite, and each member of $X_2$ belongs to all but finitely many $\bigcap\cU_n$, $n>1$.
Let $\cV_2=\cU_2$.

Step $k$: By Lemma \ref{capinfgg}, we may thin out further the families $\cU_n$, $n\ge k$, so that
they remain infinite, and each member of $X_k$ belongs to all but finitely many $\bigcap\cU_n$, $n\ge k$.
Let $\cV_k=\cU_k$.

The sets $\cV_1,\cV_2,\dots$ are as required.
\epf
Now, let $X=\Union_nX_n$ be an increasing union, where each $X_n$ satisfies
$\sone(\Omega,\Gamma)$. Let $\cU_1,\cU_2,\dots\in\Omega(X)$. As in the argument of the proof of
Lemma \ref{capinfgg}, we may assume that $\cU_n\in\Omega(X_n)$ for all $n$.
As $X_n$ satisfies $\sone(\Omega,\Gamma)$, we may thin out $\cU_n$ so that $\cU_n\in\Gamma(X_n)$.\footnote{Clearly, if
$\cU\in\Omega(X)$ and $X$ satisfies $\sone(\Omega,\Gamma)$, then $\cU$ contains a subcover $\cV$ such that
$\cV\in\Gamma(X)$ \cite{GN}.}
Apply Lemma \ref{jlem}, and pick for each $n$ some $U_n\in\cV_n\sm\{U_1,\dots,U_{n-1}\}$.
$\sseq{U_n}\in\Gamma(X)$.
\epf

Let $C(X)$ be the family of continuous functions $f:X\to\R$.
$C(X)$ has the \emph{Arhangel'ski\u{\i} property $\alpha_2$} if the following
holds for all $f^n_m\in C(X)$, $n,m\in\N$:
If for each $n$, $\lim_{m\to\oo} f^n_m(x)=0$ for all $x\in X$, then there are $m_n$ such that
$\lim_n f^n_{m_n}(x)=0$ for each $x\in X$.

For the same reason briefly mentioned in the proof of Lemma \ref{capinfgg},
one may require in the definition of $\alpha_2$ that there are infinite $I_1,I_2,\dots\sbst\N$,
such that some (equivalently, any) enumeration of the countable set $\Union_n\{f^n_m(x) : m\in I_n\}$ converges to $0$,
for each $x\in X$. Indeed, this was Arhangel'ski\u{\i}'s original definition of $\alpha_2$.

\bdfn
Let $X=\Union_nX_n$ be an increasing union. $C(X)$ is \emph{$\seq{X_n}$-$\alpha_2$} if
the following holds whenever $f^n_m\in C(X_n)$ for all $m,n$:
If for each $n$, $\lim_m f^n_m(x)=0$ for all $x\in X_n$,
then there are infinite $I_1,I_2,\dots\sbst\N$, such that
for each $x\in X$, $\lim_n\sup\{|f^n_m(x)| : m\in I_n\}=0$.\footnote{For each $x\in X$, $x\in X_n$ for all large enough $n$.
Thus, for all large enough $n$, $r_n(x)=\sup\{|f^n_m(x)| : m\in I_n\}$ is defined, and therefore the question, whether
$\lim_nr_n(x)$ converges to $0$ or not, makes sense.}
\edfn

Clearly, $\seq{X_n}$-$\alpha_2$ implies $\alpha_2$.
The proof of Lemma \ref{jlem}, with minor modifications, yields the following
new result about function spaces.

\bthm\label{jalpha}
Let $X=\Union_nX_n$ be an increasing union, such that $C(X_n)$ is $\alpha_2$ for all $n$.
Then $C(X)$ is $\seq{X_n}$-$\alpha_2$.
\ethm
\bpf
Assume that $f^n_m\in C(X_n)$ for all $m,n$, and for each $n$, $\lim_m f^n_m(x)=0$ for all $x\in X_n$.

Step 1: As $C(X_1)$ is $\alpha_2$, there are infinite $J_1,J_2,\dots\sbst\N$ such that
any enumeration of the countable set $\Union_n\{f^n_m(x) : m\in J_n\}$ converges to $0$,
for each $x\in X_1$.
In particular, $\lim_n\sup\{|f^n_m(x)| : m\in J_n\}=0$ for all $x\in X_1$.
Let $I_1=J_1$.

Step $k$: For each $n\geq k$, $\lim_{m\in J_n} f^n_m(x)=0$ for all $x\in X_k$.
As $C(X_k)$ is $\alpha_2$, the index sets $J_k,J_{k+1},\dots$ may be thinned out
so that they remain infinite, and $\lim_n\sup\{|f^n_m(x)| : m\in J_n\}=0$ for all $x\in X_k$.
Let $I_k=J_k$.

The index sets $I_1,I_2,\dots$ are as required.
\epf

\bcor
Assume that for each subset $Y$ of $X$, $C(Y)$ is $\alpha_2$.
For all $X_1\sbst X_2\sbst\dots\sbst X$, $C(\Union_nX_n)$ is $\seq{X_n}$-$\alpha_2$.\qed
\ecor

A topological space is \emph{perfectly normal} if for each
pair of disjoint closed sets $C_0,C_1\sbst X$, there is $f\in C(X)$ such that $f\inv[\{0\}]=C_0$ and $f\inv[\{1\}]=C_1$.

\bcor
Assume that $X$ is a perfectly normal space and $C(X)$ is $\alpha_2$.
For all closed sets $X_1\sbst X_2\sbst\dots\sbst X$, $C(\Union_nX_n)$ is $\seq{X_n}$-$\alpha_2$.
\ecor
\bpf
If $C(X)$ is $\alpha_2$ and $Y$ is a closed subset of $X$, then $C(Y)$ is $\alpha_2$.\footnote{The argument is
as in Theorem 4.1 of \cite{BRR91}: As $X$ is perfectly normal, the open set $X\sm Y$ is a countable increasing union of
closed sets, $X\sm Y=\Union_n C_n$. For each $n$, extend each $f^n_m$ to an element of $C(X)$
which is constantly $0$ on $C_n$. Applying $\alpha_2$ on the new sequences, we obtain $m_n$
such that $\lim_n f^n_{m_n}(x)=0$ for all $x\in X$, and thus this is also the case for
the original functions, for all $x\in Y$.}
Apply Theorem \ref{jalpha}.
\epf

In addition to these applications, Jordan's Lemma \ref{jlem} is an ingredient
in a solution of a thus far open problem, to which we now turn.

\section{$\gamma$-sets of reals from a weak hypothesis}\label{JMSSSec}

In this section, we construct sets of reals satisfying $\sone(\Omega,\Gamma)$.
Traditionally, general topological spaces satisfying $\sone(\Omega,\Gamma)$ are called
\emph{$\gamma$-spaces}, and if they happen to be (homeomorphic to) sets of real numbers,
they are called \emph{$\gamma$-sets}.

The problem settled by our construction has some history, which we now survey briefly.
This involves combinatorial cardinal characteristics of the continuum \cite{BlassHBK}.
We give the necessary definitions as we proceed.
Readers who are new to this field, may skip this section in their first reading.

$\gamma$-spaces were introduced by Gerlits and Nagy in \cite{GN},
their most influential paper, as the third property
in a list numbered $\alpha$ through $\epsilon$. This turned out to be the most important
property in the list, and obtained its alphabetic number as it name.
One of the main results in \cite{GN} is that for Tychonoff spaces $X$,
$C(X)$ with the topology of pointwise convergence is Fr\'echet-Urysohn if, and
only if, $X$ is a $\gamma$-space.

While uncountable $\gamma$-spaces exist in ZFC\footnote{The axioms of Zermelo and Fraenkel, together with the axiom of Choice,
the ordinary axioms of mathematics.} \cite{Tod95},
Borel's Conjecture (which is consistent with, but not provable within, ZFC)
implies that all \emph{metrizable} $\gamma$-spaces are countable.

Since we are dealing with constructions rather than general results,
we restrict attention in this section to subsets of $\R$
(or, since the property is preserved by continuous images,
subsets of any topological space which can be embedded in $\R$).
Thus, we may restrict attention to \emph{countable} open covers.
As mentioned above, by \emph{$\gamma$-set} we mean a $\gamma$-space which is
(homeomorphic to) a set of real numbers.

Gerlits and Nagy proved in \cite{GN} that
Martin's Axiom implies that all sets of reals of cardinality less than $\fc$ are $\gamma$-sets.
There is a simple reason for that: The \emph{critical cardinality} of a property $P$,
denoted $\non(P)$, is the minimal cardinality of a set not satisfying $P$.
Let $\binom{\Omega}{\Gamma}$ be the property: Each $\cU\in\Omega(X)$ contains a
set $\cV\in\Gamma(X)$. Gerlits and Nagy proved that $\sone(\Omega,\Gamma)=\binom{\Omega}{\Gamma}$ \cite{GN}.
Let $A\as B$ mean that $A\sm B$ is finite.
$A$ is a \emph{pseudointersection} of $\cF$ if $A\as B$ for all $B\in\cF$.
Let $\fp$ be the minimal cardinality of a family
$\cF$ of infinite subsets of $\N$ which is closed under finite intersections,
and has no pseudointersection.
Then $\non\binom{\Omega}{\Gamma}=\fp$ \cite{GM}, and Martin's Axiom implies $\fp=\fc$ \cite{GM}.

By definition, for each property $P$ of sets of reals, every set of reals whose cardinality is
smaller than $\non(P)$ satisfies $P$. Thus, the real question is whether there is a set of reals $X$
of cardinality at least $\non(P)$, which satisfies $P$.
Galvin and Miller \cite{GM} proved a result of this type: $\fp=\fc$
implies that there is a $\gamma$-set of cardinality $\fp$.
Just, Miller, Scheepers and Szeptycki \cite{coc2} have improved the construction
of \cite{GM}.
We introduce their construction in a slightly more general form, that will be useful later.

Cantor's space $\Cantor$ is equipped with the Tychonoff product topology,
and $P(\N)$ is identified with $\Cantor$ using characteristic functions.
This defines the topology of $P(\N)$. The partition $P(\N)=\roth\cup\Fin$, into
the infinite and the finite sets, respectively, is useful here.

For $f,g\in\NN$, let $f\le^* g$ if $f(n)\le g(n)$ for all but finitely many
$n$. $\fb$ is the minimal cardinality of a $\le^*$-unbounded subset of $\NN$.
A set $B\sbst\roth$ is \emph{unbounded} if the set of all increasing enumerations
of elements of $B$ is unbounded in $\NN$, with respect to $\le^*$.
It follows that $|B|\ge\fb$.
For $m,n\in\N$, let $(m,n)=\{k : m<k<n\}$.

\blem[folklore]\label{infdisj}
If $B\sbst\roth$ is unbounded, then for each increasing $f\in\NN$,
there is $x\in B$ such that $x\cap (f(n),f(n+1))=\emptyset$ for infinitely
many $n$.
\elem
\bpf
Assume that $f$ is a counterexample. Let $g$ dominate all functions $f_m(n)=f(n+m)$, $m\in\N$.
Then for each $x\in B$, $x\le^* g$. Indeed, let $m$ be such that for all $n\ge m$,
$x\cap (f(n),f(n+1))\neq\emptyset$. Then for each $n$, the $n$-th element of $x$
is smaller than $f_{m+1}(n)$.
\epf

\bdfn
A\emph{tower} of cardinality $\kappa$ is a set $T\sbst\roth$ which can be enumerated
bijectively as $\{x_\alpha : \alpha<\kappa\}$,
such that for all $\alpha<\beta<\kappa$, $x_\beta\as x_\alpha$.

An \emph{unbounded tower} of cardinality $\kappa$ is an unbounded set $T\sbst\roth$ which is
a tower of cardinality $\kappa$. (Necessarily, $\kappa\ge\fb$.)
\edfn

Let $\ft$ be the minimal cardinality of a tower which has no pseudointersection.
Rothberger proved that $\ft\le\fb$ \cite{BlassHBK}.

\blem[folklore]\label{ut}
$\ft=\fb$ if, and only if, there is an unbounded tower of cardinality $\ft$.
\elem
\bpf
$(\Impl)$ Construct $x_\alpha$ by induction on $\alpha$.
Let $\{b_\alpha : \alpha<\fb\}\sbst\NN$ be unbounded.
At step $\alpha$, let $a$ be a pseudointersection
of $\{x_\beta : \beta<\alpha\}$, and take $x_\alpha\sbst a$
such that the increasing enumeration of $x_\alpha$ dominates $b_\alpha$.

$(\Leftarrow)$ $\ft\le\fb\le|T|=\ft$.
\epf

Just, Miller, Scheepers and Szeptycki \cite{coc2} proved that if
$T$ is an unbounded tower of cardinality $\aleph_1$,
then $T\cup\Fin$ satisfies $\sone(\Omega,\Omega)$, as well as
a property, which was later proved by Scheepers \cite{wqn} to be equivalent to $\sone(\Gamma,\Gamma)$ .
In Problem 7 of \cite{coc2}, we are asked the following.

\bprb[Just-Miller-Scheepers-Szeptycki \cite{coc2}]\label{JMSSProb}
Assume that $T\sbst\roth$ is an unbounded tower of cardinality $\aleph_1$ (so that $\aleph_1=\fb$).
Is $T\cup\Fin$ a $\gamma$-set, i.e., satisfies $\sone(\Omega,\Gamma)$?
\eprb

Scheepers proves in \cite{alpha_i} that for each unbounded tower $T$ of cardinality $\ft=\fb$,
$T\sbst\roth$ satisfies $\sone(\Gamma,\Gamma)$.

Miller \cite{MillerBC} proves that in
the Hechler model, there are no uncountable $\gamma$-sets. In this model, $\aleph_1=\fp=\ft<\fb$,
and thus $\aleph_1=\ft$ does not suffice to have an uncountable $\gamma$-set.
At the end of \cite{MillerBC} and in its appendix, Miller proves that $\lozenge(\fb)$,
a property strictly stronger than $\aleph_1=\fb$, implies that there is an uncountable $\gamma$-set.\footnote{$\lozenge(\fb)$
is defined in Dzamonja-Hrusak-Moore \cite{DHM}.}
He concludes that it is still open whether $\fb = \aleph_1$ is enough to construct an uncountable $\gamma$-set.

We show that the answer is positive, and indeed also answer a question of Gruenhage and Szeptycki \cite{FUfin}:
A classical problem of Malykhin asks whether there is a countable Fr\'echet-Urysohn topological group which is not
metrizable. Gruenhage and Szeptycki prove that $F\sbst\NN$ is a $\gamma$-set if, and only if, a certain construction
associated to $F$ provides a positive answer to Malykhin's Problem \cite{FUfin}.
They define a generalization of $\gamma$-set, called \emph{weak $\gamma$-set}, and combine their results with results of Nyikos
to prove the following.

\bthm[Gruenhage-Szeptycki \cite{FUfin}]
If $\fp=\fb$, then there is a weak $\gamma$-set in $\NN$.
\ethm

They write: ``The relationship between $\gamma$-sets and weak $\gamma$-sets is not known.
Perhaps $\fb = \fp$ implies the existence of a $\gamma$-set'' \cite{FUfin}.
Our solution confirms their conjecture.

$\fp\le\ft\le\fb$, and in all known models of set theory, $\fp=\ft$.
When $\fp=\fc$, our theorem shows that in Galvin and Miller's construction from \cite{GM},
even if the possible open covers are not considered at all,
the resulting set is still a $\gamma$-set.

\bthm\label{JMSSSol}
For each unbounded tower $T$ of cardinality $\fp$ in $\roth$,
$T\cup\Fin$ satisfies $\sone(\Omega,\Gamma)$.
\ethm
\bpf
By Lemma \ref{ut}, there is an unbounded tower $T$ of cardinality $\fp$ if, and only if, $\fp=\fb$.
Let $T=\{x_\alpha : \alpha<\fb\}$ be an unbounded tower of cardinality $\fb$.
For each $\alpha$, let $X_\alpha=\{x_\beta : \beta<\alpha\}\cup\Fin$.
We will show that $T\cup\Fin$ satisfies $\binom{\Omega}{\Gamma}$. Let $\cU\in\Omega(T\cup\Fin)$.
We use the following modification of Lemma 1.2 of \cite{GM}.

\blem\label{GM+}
Assume that $\Fin\sbst X\sbst P(\N)$ and $X$ satisfies $\binom{\Omega}{\Gamma}$.
For each $\cU\in\Omega(X)$, there are $m_1<m_2<\dots$ and distinct $U_1,U_2,\dots\in\cU$
such that $\sseq{U_n}\in\Gamma(X)$, and for each $x\sbst\N$, $x\in U_n$ whenever $x\cap (m_n,m_{n+1})=\emptyset$.
\elem
\bpf
As $X$ satisfies $\binom{\Omega}{\Gamma}$, we may thin out $\cU$ so that $\cU\in\Gamma(X)$.

We proceed as in the proof of Lemma 1.2 of \cite{GM}.
Let $m_1=1$. For each $n\ge 1$: As $\cU\in\Omega(X)$, each finite subset of $X$ is contained in
infinitely many elements of $\cU$. Take $U_n\in\cU\sm\{U_1,\dots,U_{n-1}\}$,
such that $P(\{1,\dots,\allowbreak m_n\})\sbst U_n$.
As $U_n$ is open, for each $s\sbst \{1,\dots,m_n\}$ there is $k_s$ such that
for each $x\in P(\N)$ with $x\cap\{1,\dots,k_s-1\}=s$, $x\in U_n$. Let $m_{n+1}=\max\{k_s : s\sbst\{1,\dots,m_n\}\}$.

As $\sseq{U_n}$ is an infinite subset of $\cU\in\Gamma(X)$, $\sseq{U_n}\in\Gamma(X)$, too.
\epf

We may assume that $\cU$ is countable.
As $\fb$ is regular, there is $\alpha_1<\fb$ such that $X_{\alpha_1}$ is not contained in any member of $\cU$.
This guarantees that $\cU\in\Omega(X_\alpha)$ for all $\alpha\ge\alpha_1$.

As $|X_{\alpha_1}|<\fp$, $X_{\alpha_1}$ satisfies $\binom{\Omega}{\Gamma}$.
As $\cU\in\Omega(X_{\alpha_1})$, by Lemma \ref{GM+} there are
$m^1_1<m^1_2<\dots$ and distinct $U^1_1,U^1_2,\dots\in\cU$
such that $\sseq{U^1_n}\in\Gamma(X_{\alpha_1})$, and for each
$x\in P(\N)$, $x\in U^1_n$ whenever $x\cap (m^1_n,m^1_{n+1})=\emptyset$.
Let $D_1=\N$.

As $\alpha_1<\fb$, $\{x_\alpha : \alpha_1<\alpha<\fb\}$ is unbounded.
By Lemma \ref{infdisj}, there is $\alpha_2>\alpha_1$ such that
$D_2=\{n : x_{\alpha_2}\cap (m^1_n,m^1_{n+1})=\emptyset\}$ is infinite.
As $|X_{\alpha_2}|<\fp$, $X_{\alpha_2}$ satisfies $\binom{\Omega}{\Gamma}$.
As $\cU\in\Omega(X_{\alpha_2})$, by Lemma \ref{GM+} there are
$m^2_1<m^2_2<\dots$ and distinct $U^2_1,U^2_2,\dots\in\cU$
such that $\sseq{U^2_n}\in\Gamma(X_{\alpha_2})$, and for each
$x\in P(\N)$, $x\in U^2_n$ whenever $x\cap (m^2_n,m^2_{n+1})=\emptyset$.
As $D_2$ is infinite, $\{U^2_n : n\in D_2\}\in\Gamma(X_{\alpha_2})$.

Continue in the same manner, to define for each $k>1$ elements with the following properties:
\be
\itm $\alpha_k>\alpha_{k-1}$;
\itm $D_k=\{n : x_{\alpha_k}\cap (m^{k-1}_n,m^{k-1}_{n+1})=\emptyset\}$ is infinite;
\itm $m^k_1<m^k_2<\dots$;
\itm $U^k_1,U^k_2,\dots\in\cU$ are distinct;
\itm $\{U^k_n : n\in D_k\}\in\Gamma(X_{\alpha_k})$; and
\itm For each $x\in P(\N)$, $x\in U^k_n$ whenever $x\cap (m^k_n,m^k_{n+1})=\emptyset$.
\ee
Let $\alpha=\sup_k\alpha_k$. As $\fb$ is regular, $\alpha<\fb$.
$X_\alpha=\Union_k X_{\alpha_k}$ is a countable increasing union.
For each $k$, $|X_{\alpha_k}|<\fb$, and thus $X_{\alpha_k}$ satisfies $\sone(\Gamma,\Gamma)$.
By Lemma \ref{jlem}, there are infinite $I_1\sbst D_1,I_2\sbst D_2,\dots$ such that
each $x\in X_\alpha$ belongs to  $\bigcap_{n\in I_k}U^k_n$ for all but finitely many $k\in\N$.

Take $n_1\in I_2$.
For $k>1$, take $n_k\in I_{k+1}$ such that $m^k_{n_k}>m^{k-1}_{n_{k-1}+1}$,
$x_\alpha\cap (m^k_{n_k},m^k_{n_k+1})\sbst x_{\alpha_{k+1}}\cap (m^k_{n_k},m^k_{n_k+1})$,
and $U^k_{n_k}\nin\{U^1_{n_1},\dots,U^{k-1}_{n_{k-1}}\}$.
We claim that $\{U^k_{n_k} : k\in\N\}\in\Gamma(T\cup\Fin)$.
As $\{U^k_{n_k} : k\in\N\}\in\Gamma(X_\alpha)$, it remains to show that for each
$x\as x_\alpha$, $x\in U^k_{n_k}$ for all but finitely many $k$.
For each large enough $k$, $m^{k}_{n_k}$ is large enough, so that
$$x\cap (m^k_{n_k},m^k_{n_k+1})\sbst x_\alpha\cap (m^k_{n_k},m^k_{n_k+1})\sbst
x_{\alpha_{k+1}}\cap (m^k_{n_k},m^k_{n_k+1})=\emptyset,$$
since $n_k\in D_{k+1}$. Thus, $x\in U^k_{n_k}$.
\epf

\brem
Zdomskyy points out that our proof actually shows that a wider family of sets are $\gamma$-sets.
For example, if we start with $T$ an unbounded tower of cardinality $\fp$, and thin out its elements
arbitrarily, $T\cup\Fin$ remains a $\gamma$-set. This may be useful for constructions of examples
with additional properties, since this way, each element of $T$ may be chosen arbitrarily from a certain perfect set.
\erem

In particular, we have that in each model of ZFC where $\fp=\fb$, there are $\gamma$-sets of cardinality $\fp$.
In the following corollary, by ``$X$ model'' we mean the model obtained by generically extending a model of \CH{}
by adding $\aleph_2$ (or more, in the case $X\in\{\mbox{Cohen, Random}\}$) $X$ reals in the standard way --
see \cite{BlassHBK} for details.

\bcor\label{ec}
In each of the Cohen, Random, Sacks, and Miller models of ZFC, there are $\gamma$-sets of reals
with cardinality $\fp$.\qed
\ecor

Earlier, Corollary \ref{ec} was shown for the Sacks model by Ciesielski, Mill\'an, and Pawlikowski
in \cite{CPAgamma},
and for the Cohen and Miller models by Miller \cite{MillerBC}, using specialized
arguments. It seems that the result, that there are uncountable $\gamma$-sets
in the Random reals model (constructed by extending a model of \CH{}), is new.
We point out that the Random poset alone is not the reason for having uncountable $\gamma$-sets
in the generic extension: Judah, Shelah, and Woodin prove in \cite{JSW} that
there are no uncountable strong measure zero sets (and in particular, no uncountable $\gamma$-sets)
in an extension of Laver's model by random reals.

As discussed above, there are no uncountable $\gamma$-sets in the Hechler model \cite{MillerBC}.
Since the Laver and Mathias models satisfy Borel's Conjecture, there are no
uncountable $\gamma$-sets in these models, too.

\section{Heredity}\label{heredsec}

A topological space $X$ satisfies a property $P$ \emph{hereditarily} if
each subspace of $X$ satisfies $P$.
In our context, heredity was observed to be tightly connected to the following
property. $X$ is a \emph{$\sigma$-space} if each Borel subset of $X$ is an $F_\sigma$
subset of $X$.

A combination of results of Fremlin-Miller \cite{FM88},
Bukovsk\'y-Rec\l{}aw-Repick\'y \cite{BRR01}, and Bukovsk\'y-Hale\v{s} \cite{BH03}
implies that a Tychonoff space
$X$ satisfies $\ufin(\rmO,\Gamma)$ hereditarily if, and only if,
$X$ satisfies $\ufin(\rmO,\Gamma)$ and is a $\sigma$-space (see also \cite{hH}).\footnote{The assumptions
on $X$ in the cited references are stronger, but: A Tychonoff hereditarily-$\ufin(\rmO,\Gamma)$ space
cannot have the unit interval $[0,1]$ as a continuous image, and thus is zero-dimensional. As it is hereditarily Lindel\"of,
each open set is a countable union of clopen sets, and this suffices for the arguments
in the cited references.}
A similar result was proved for $\sone(\Gamma,\Gamma)$ in \cite{BH03}.
Problem 7.9 in \cite{BRR01} asks whether every (nice enough)
$\sigma$-space $X$ satisfying $\sone(\Omega,\Gamma)$, satisfies
$\sone(\Omega,\Gamma)$ hereditarily.
A negative answer was given by Miller \cite{MilNonGamma}.

We have, in light of the previous section, a simple reason for the difference
between $\sone(\Omega,\Gamma)$ on one hand, and $\sone(\Gamma,\Gamma)$ and $\ufin(\rmO,\Gamma)$
on the other hand.

\bdfn
Let $\scrA\in\{\Gamma,\Omega,\rmO,\dots\}$. $\scrA$ is \emph{Borel superset covering} for $X$
if, for each subspace $Y\sbst X$ and each $\cU\in\scrA(Y)$, there are $\cV\sbst\cU$ and a Borel $B\sbst X$,
such that $Y\sbst B$ and $\cV\in\scrA(B)$.
\edfn

Many classical types of covers are Borel superset covering. In particular, we have the following.

\blem\label{example}
$\Gamma$ and $\rmO$ are Borel superset covering (for all $X$).
\elem
\bpf
If $\cU\in\Gamma(X)$, take a countable infinite subset $\cV=\sseq{U_n}\sbst\cU$, and $B=\Union_m\bigcap_{n\ge m}U_n$.
If $\cU\in\rmO(X)$, take $B=\Union\cU$.
\epf

By Miller's mentioned result and the following theorem, $\Omega$ need not be Borel superset covering.

\bthm\label{ssCtoH}
Let $\scrA,\scrB\in\{\Gamma,\Omega,\rmO\}$, and $\Pi\in\{\sone,\sfin,\ufin\}$.
Assume that $X$ is a $\sigma$-space, satisfies $\Pi(\scrA,\scrB)$, and
$\scrA$ is Borel superset covering for $X$.
Then $X$ satisfies $\Pi(\scrA,\scrB)$ hereditarily.
\ethm
\bpf
Let $Y\sbst X$ and assume that $\cU_1,\cU_2,\dots\in\scrA(Y)$. For each $n$, pick intermediate
Borel sets $Y\sbst B_n\sbst X$ such that $\cU_n\in\scrA(B_n)$, and let $B=\bigcap_nB_n$.
Then $\cU_1,\cU_2,\dots\in\scrA(B)$.
As $B$ is Borel and $X$ is a $\sigma$-space, $B$ is $F_\sigma$ in $X$. By Corollary \ref{Fs},
$B$ satisfies $\Pi(\scrA,\scrB)$, and we can thus obtain from the covers $\cU_1,\cU_2,\dots$
a cover $\cV\in\scrB(B)$. Then $\cV\in\scrB(Y)$.
\epf

As for the other direction, we have, by the discussion at the beginning of this section,
the following.

\blem\label{crit}
If a property $P$ of Tychonoff spaces implies $\ufin(\rmO,\Gamma)$
and $X$ satisfies $P$ hereditarily, then $X$ is a $\sigma$-space.\qed
\elem

During his work with the second author on \cite{hH}, L. Zdomskyy asked
whether every subset of $\R$ satisfying $\ufin(\rmO,\Omega)$ hereditarily
is a $\sigma$-space.
We show that the answer is ``No'', in a very strong sense.
We will see that, in the context of Scheepers Diagram, Lemma \ref{crit}
becomes a \emph{criterion}, that is, satisfying $P$ hereditarily implies
being a $\sigma$-space if, and only if, $P$ implies $\ufin(\rmO,\Gamma)$.

We use, in our proof, the method of ``forcing'', following an elegant approach
of Brendle \cite{Brendle96}. A proof avoiding this method (e.g., using the
methods of \cite{o-bdd}) would be more lengthy and technically involved.

\bthm\label{OmOm}
Assume \CH. There is $L\sbst\R$ such that all finite powers of $L$
satisfy $\sone(\BO,\BO)$ hereditarily, but $L$ is not a $\sigma$-space.
\ethm
\bpf
Let $M_\alpha$, $\alpha < \aleph_1$, be an increasing
sequence of countable submodels of $H(\lambda)$ ($\lambda$ large enough)
such that:
\be
\itm $\R\sbst \bigcup_{\alpha<\aleph_1} M_\alpha$;
\itm For each $\alpha < \aleph_1$, $M_{\alpha+1} \models M_\alpha$ is countable; and
\item $(M_\beta: \beta\leq \alpha) \in M_{\alpha+1}$.
\ee
For each $\alpha$ choose a real $c_\alpha$ Cohen generic over
$M_\alpha$. (As $M_\alpha$ is countable, there is such $c_\alpha$.)
Let $L=\{c_\alpha: \alpha <\aleph_1\}$.
It is well known that a set $L$ constructed this way is a \emph{Luzin set}, i.e.,
$L$ is uncountable, but for each meager $M\sbst\R$, $L\cap M$ is countable.\footnote{Indeed,
each meager set is contained in a meager $F_\sigma$ set $M$, which in turn is coded by a single
real $r$. Let $\alpha$ be such that $r\in M_\alpha$. Then for each $\beta>\alpha$,
$c_\beta\nin M$.}

\blem[folklore]\label{LuzinNotSigma}
No Luzin set $L$ is a $\sigma$-space.
\elem
\bpf
Assume otherwise, and take a countable dense $D\sbst L$. As $L\sm D$ is Borel, it is $F_\sigma$.
Then $D=\bigcap_nU_n$, an intersection of open sets. Then for each $n$, $L\sm U_n$ is meager,
and thus countable. It follows that $L\sm D$ is countable, a contradiction.
\epf

Recall that according to Definition \ref{pows}, $X$ satisfies $P^\uparrow$ if
all finite powers of $X$ satisfy $P$.
By \cite{ideals}, or alternatively Lemma \ref{example} and Remark \ref{BorHered},
$\sone(\B,\B)$ is hereditary.
By \cite{CBC}, $\sone(\BO,\BO)^\uparrow=\sone(\B,\B)^\uparrow$.

\blem\label{powhered}
If $P$ is hereditary, then so is $P^\uparrow$.
\elem
\bpf
If $Y\sbst X$ and $X$ satisfies $P^\uparrow$, then
for each $k$, $Y^k\sbst X^k$ and $X^k$ satisfies $P$.
Thus, $Y^k$ satisfies $P$ for all $k$.
\epf

Thus, it suffices to prove that all finite powers of $L$ satisfy $\sone(\B,\B)$.
Consider for example $L^2$. We have to show (essentially, \cite{CBC}) that for each
Borel $\Psi:\R^2\to\NN$, there is $g\in\NN$ such that for each $(x,y)\in L^2$,
$\Psi(x,y)(n)=g(n)$ for some $n$.

Let $\Psi: \R^2\to\NN$ be Borel.
$\Psi$ is coded by a real. Let $\alpha$ be such that this code belongs to $M_{\alpha}$.
$\NN\cap M_\alpha$ is countable. Take a canonical partition $\N=\Union_nI_n$ with each
$I_n$ infinite, so that $(I_n : n\in\N)\in M_\alpha$. Enumerate $\NN\cap M_\alpha=\{f_n : n\in\N\}$, and take $g\in\NN$ such that
for each $n$, $g\mid_{I_n} = f_n\mid_{I_n}$.
For $x,y\in L$, $\Psi(x,y)\in M_\alpha[x,y]$, an extension of $M_\alpha$ by finitely many Cohen reals,
which is in fact either $M_\alpha$ or an extension of $M_\alpha$ by a single Cohen real.
As the Cohen forcing does not add eventually different reals (e.g., \cite{barju}), there is $n$ such that $f_n\mid_{I_n}$
coincides with $\Psi(x,y)\mid_{I_n}$ on infinitely many values. But $f_n\mid_{I_n}=g\mid_{I_n}$.
\epf

We conclude the section with the following characterization, which may be viewed
as a necessary revision of Problem of 7.9 of Bukovsk\'y-Rec\l{}aw-Repick\'y \cite{BRR01},
so that it has a (provably) positive answer.\footnote{See the discussion at the beginning
of the present section.}

\bthm\label{brr}
A Tychonoff space $X$ satisfies $\sone(\Omega,\Gamma)$ hereditarily if, and only if:
\be
\itm $X$ satisfies $\sone(\Omega,\Gamma)$;
\itm $X$ is a $\sigma$-space; and
\itm $\Omega$ is Borel superset covering for $X$.
\ee
\ethm
\bpf
Theorem \ref{ssCtoH} provides the ``if'' part.
For the ``only if'' part, by Lemma \ref{crit} it remains to prove (3).
Let $Y\sbst X$, $\cU\in\Omega(Y)$. As $Y$ satisfies $\sone(\Omega,\Gamma)$, there are
$U_n\in\cU$ such that $\sseq{U_n}\in\Gamma(Y)$. Let $B=\Union_m\bigcap_{n\ge m}U_n$.
Then $Y\sbst B$ and $\sseq{U_n}\in\Omega(B)$.
\epf

\brem[The case of countable Borel covers]\label{BorHered}
In the Borel case, the properties are hereditary for Borel subsets \cite{CBC},
and thus in Theorem \ref{ssCtoH}, there is no need to assume that $X$ is a $\sigma$-space.
This explains why some of the Borel properties are hereditary \cite{ideals} whereas others are not \cite{MilNonGamma},
and why none of the open properties is hereditary \cite{ideals}.

Two major open problems concerning $\sone(\BO,\BO)$ and $\sfin(\BO,\BO)$ are whether these properties
are hereditary, and whether they are---like their open variant---preserved by finite powers \cite{ideals, OPiT, MilNonGamma}.
By Lemma \ref{powhered}, the problems are related: If $\sone(\BO,\BO)$ is preserved by finite powers,
then it is equal to $\sone(\B,\B)^\uparrow$, which is hereditary. And similarly for $\sfin(\BO,\BO)$.

Item (2) may be removed from
Theorem \ref{brr} in the Borel case. On the other hand, $\ufin(\B,\BG)$ is hereditary and (for Lindel\"of spaces)
implies $\ufin(\rmO,\Gamma)$, and thus implies being a $\sigma$-set. Theorem \ref{OmOm} completes the picture.
\erem

\section{Topological groups with strong combinatorial properties}\label{topgpssec}

Problem 10.7 in \cite{OPiT}, which is also implicit in Theorem 20 of \cite{o-bdd} and in the discussion
around it, asks whether \CH{} implies the existence of an uncountable subgroup $G$ of $\Bgp$
satisfying $\sone(\Omega,\Gamma)$. We will show that the answer is ``Yes'', and indeed
the weaker hypothesis $\fp=\fb$ suffices to obtain such groups (see the discussion in Section \ref{JMSSSec}).
Indeed, Theorem \ref{JMSSSol} gives an uncountable subset of $\Bgp$ satisfying $\sone(\Omega,\Gamma)$.
We will show that for a wide class of properties $P$,
including $\sone(\Omega,\Gamma)$, if $X$ satisfies $P$, then so is the group generated
by $X$ in any Tychonoff topological group $G$.

There is some restriction on $P$:
Let $F(X)$ be the free topological group generated by $X$.
That is, the group such that each continuous function $f$ from $X$ into a topological group $H$,
can be extended uniquely to a continuous homomorphism from $F(X)$ to $H$.\footnote{A thorough introduction
to free topological groups is available in \cite{TkaTGII}.}

\blem
Assume that $P$ is a property of Tychonoff topological spaces, which is hereditary for closed subsets
and preserved by homeomorphic images. If for each $X$ satisfying $P$, $F(X)$ satisfies $P$, then $P=P^\uparrow$.
\elem
\bpf
Assume that $X$ satisfies $P$. Then $F(X)$ satisfies $P$.
For each $k$, $X^k$ is homeomorphic to $\{x_1\cdots x_k : x_1,\dots,x_k\in X\}$, a closed subspace of $F(X)$, and thus satisfies $P$.
\epf

Thus, we must consider properties $P$ such that $P=P^\uparrow$.
For each $Q$ in the Scheepers Diagram, $P=Q^\uparrow$ is as required.
By the discussion in the introduction \cite{coc2}:
\be
\itm $\sone(\rmO,\rmO)^\uparrow = \sone(\Omega,\Omega)^\uparrow = \sone(\Omega,\Omega)$;
\itm $\sfin(\rmO,\rmO)^\uparrow = \ufin(\rmO,\Omega)^\uparrow = \sfin(\Gamma,\Omega)^\uparrow = \sfin(\Omega,\Omega)^\uparrow = \sfin(\Omega,\Omega)$.
\ee
Moreover, $\sone(\Omega,\Gamma)^\uparrow = \sone(\Omega,\Gamma)$ \cite{GM}, and $\ufin(\rmO,\Gamma)^\uparrow$ also has a simple
characterization \cite{coc7}. The properties $\sone(\Gamma,\Gamma)^\uparrow$, $\sone(\Gamma,\Omega)^\uparrow$, and $\sone(\Gamma,\rmO)^\uparrow$
seem, however, to be unexplored. By the results of \cite{coc2}, each of these properties is strictly stronger than its non-$\uparrow$-ed version.

All $P$ in the Scheepers Diagram have the properties required in the following theorem.

\bthm\label{gengp}
Assume that $P$ is a property of Tychonoff topological spaces, which is hereditary for closed subsets and preserved by continuous images.
If $P$ is linearly $\sigma$-additive, then for each $X$ satisfying $P^\uparrow$
and each Tychonoff topological group $G$ containing $X$,
the group $\<X\>\le G$ generated by $X$ satisfies $P^\uparrow$.
\ethm
\bpf
It is not difficult to prove the following.
\blem
Assume that $P$ is a property of topological spaces, which is hereditary for closed subsets and preserved by continuous images,
and $P$ is linearly $\sigma$-additive. Then $P^\uparrow$ also has these three properties.\qed
\elem

Fix distinct $a,b\in X\sm\{1\}$, $1$ being the identity element of $G$.
As $\{1,a\},\{a,b\}$ are discrete subspaces of $G$, $X\x\{1,a\}$ is homeomorphic to $X\x\{a,b\}$.
As $X^2$ satisfies $P^\uparrow$, so does its closed subset $X\x\{a,b\}$, and thus
so does $X\x\{1,a\}$. Thus, so does the image $Y$ of $X\x\{1,a\}$ under the continuous map
$(x,y)\mapsto xy\inv$. $Y=X\cup Xa\inv$, and as $a\in X$, $1\in Y$.
As $a\in X$, $\<X\>=\<Y\>$. Thus, we may assume that $1\in X$.

$(\{1\}\x X)\cup (X\x\{1\})$ is a closed subset of $X^2$,
which satisfies $P^\uparrow$. Thus, $(\{1\}\x X)\cup (X\x\{1\})$ satisfies $P^\uparrow$,
and therefore so does its image under the same map, $Y=X\cup X\inv$.
As $\<Y\>=\<X\>$, we may assume that $X=X\inv$, and $1\in X$.

Thus,
$$\langle X\rangle=\Union_n \{x_1\cdots x_n : x_1,\dots,x_n\in X\}$$
is an increasing union. For each $n$, $\{x_1\cdots x_n : x_1,\dots,x_n\in X\}$ is a continuous
image of $X^n$, which satisfies $P^\uparrow$. Thus,
$\<X\>$ satisfies $P^\uparrow$.
\epf

Theorem \ref{gengp} can also be stated in the language of free topological groups.
We obtain results analogous (but incomparable) to ones of Banakh, Repov\v{s}, and Zdomskyy \cite{FrGps}.

\bthm\label{frgengp}
Let $P$ be as in Theorem \ref{gengp}.
For each topological space $X$ satisfying $P^\uparrow$, the
free topological group $F(X)$ satisfies $P^\uparrow$.\qed
\ethm

It follows that for each $P$ in the Scheepers Diagram, and each $X$ satisfying $P^\uparrow$,
$\<X\>$ satisfies $P^\uparrow$. Previous constructions
of topological groups satisfying $P^\uparrow$ for these properties were much more involved than constructions of topological
\emph{spaces} satisfying $P^\uparrow$.

By Theorem \ref{JMSSSol} and Theorem \ref{gengp}, we have the following.

\bcor
Assume that $\fp=\fb$. There is a subgroup of $\R$ of cardinality $\fb$, satisfying $\sone(\Omega,\Gamma)$.\qed
\ecor

\begin{figure}[!htp]
\begin{changemargin}{-3cm}{-3cm}
\begin{center}
{%\tiny
%\scriptsize
$\xymatrix@C=-2pt@R=8pt{%@=7pt{
%1
&
&
& \nsr{\ufin(\rmO,\Gamma)}{\fb}\ar[rr]%\ar@{.>}[dr]^?
&
& \nsr{\ufin(\rmO,\Omega)}{\fd}\ar[rrrrr]%\ar@/_/@{.>}[dl]_?
&
&
&
&
&
&
& \nsr{\sfin(\rmO,\rmO)}{\fd}
\\
%2
&
&
&
& \nsr{\sfin(\Gamma,\Omega)}{\fd}\ar[ur]
\\
%3
& \nsr{\sone(\Gamma,\Gamma)}{\fb}\ar[rr]\ar[uurr]
&
& \nsr{\sone(\Gamma,\Omega)}{\fd}\ar[rrr]\ar[ur]
&
&
& \nsr{\sone(\Gamma,\rmO)}{\fd}\ar[uurrrrrr]
\\
%4
  \nsr{\sone(\BG,\BG)}{\fb}\ar[ur]\ar[rr]
&
& \nsr{\sone(\BG,\BO)}{\fd}\ar[ur]\ar[rrr]
&
&
& \nsr{\sfin(\B,\B)}{\fd}\ar[ur]
\\
%5
&
&
&
& \nsr{\sfin(\Omega,\Omega)}{\fd}\ar'[u]'[uu][uuu]
\\
%6
\\
%7
&
& \nsr{\sfin(\BO,\BO)}{\fd}\ar[uuu]\ar[uurr]
\\
%8
& \nsr{\sone(\Omega,\Gamma)}{\fp}\ar'[r][rr]\ar'[uuuu][uuuuu]
&
& \nsr{\sone(\Omega,\Omega)}{\cov(\cM)}\ar'[uuuu][uuuuu]\ar'[rr][rrr]\ar[uuur]
&
&
& \nsr{\sone(\rmO,\rmO)}{\cov(\cM)}\ar[uuuuu]
\\
%9
  \nsr{\sone(\BO,\BG)}{\fp}\ar[uuuuu]\ar[rr]\ar[ur]
&
& \nsr{\sone(\BO,\BO)}{\cov(\cM)}\ar[uu]\ar[ur]\ar[rrr]
&
&
& \nsr{\sone(\B,\B)}{\cov(\cM)}\ar[uuuuu]\ar[ur]
}$
}
\caption{The extended Scheepers Diagram}\label{extSch}
\end{center}
\end{changemargin}
\end{figure}
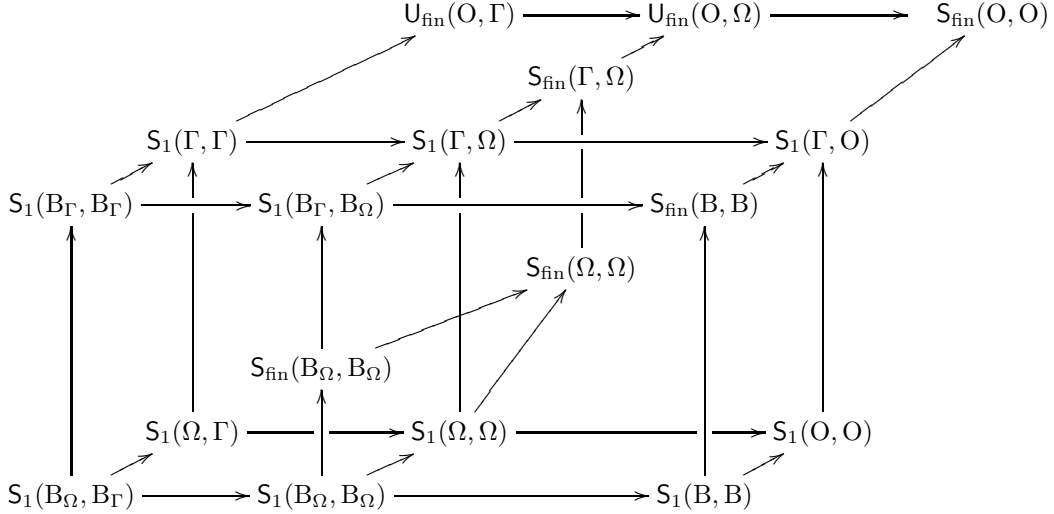

We conclude with an application of the results of this paper.
To put it in context, we draw in Figure \ref{extSch} the Scheepers Diagram, extended to also contain the
Borel properties \cite{CBC}. In the Borel case, some additional equivalences hold, and thus there
are fewer distinct properties. In particular, $\ufin(\B,\BG)$, the Borel counterpart of $\ufin(\rmO,\Gamma)$,
is equivalent to $\sone(\BG,\BG)$ \cite{CBC}.

\bthm\label{last}
Assume \CH{}.
\be
\itm There is a subgroup $G$ of $\R$ of cardinality continuum,
such that each finite power of $G$ satisfies $\sone(\BO,\BO)$ hereditarily, but
$G$ is not a $\sigma$-space, and does not satisfy $\ufin(\rmO,\Gamma)$.
\itm There is a subgroup $H$ of $\R$ of cardinality continuum,
such that all finite powers of $H$ satisfy $\sone(\BG,\BG)$ and $\sfin(\BO,\BO)$ hereditarily, but $H$ does not satisfy $\sone(\rmO,\rmO)$.
\ee
\ethm
\bpf
By \cite{ideals}, or alternatively Lemma \ref{example} and Remark \ref{BorHered},
$\sone(\B,\B)$, $\sone(\BG,\BG)$, and $\sfin(\B,\B)$, are all hereditary.
By \cite{CBC}, $\sone(\BO,\BO)^\uparrow=\sone(\B,\B)^\uparrow$, and similarly for $\sfin$.
Thus, by Lemma \ref{powhered}, it suffices to prove the statements with the words ``hereditarily'' removed.

(1) In the proof of Theorem \ref{OmOm}, we have constructed a Luzin set $L\sbst\R$ satisfying $\sone(\BO,\BO)^\uparrow$.
$L$ is not a $\sigma$-space (Lemma \ref{LuzinNotSigma}).
Let $G=\<L\>$.
By Theorem \ref{gengp}, $G$ satisfies $\sone(\BO,\BO)^\uparrow$.

As $L\sbst G$ and being a $\sigma$-space is hereditary, $G$ is not a $\sigma$-space.
As $\sone(\B,\B)$ is hereditary and is not satisfied by perfect subsets of $\R$,\footnote{Assume
that a perfect set satisfies $\sone(\B,\B)$. Then so does the unit interval $[0,1]$, its continuous image.
But even $\sone(\rmO,\rmO)$ implies Lebesgue measure zero.}
$G$ does not contain any perfect subset of $\R$.
Assume that $G$ satisfies $\ufin(\rmO,\Gamma)$. Theorem 5.5 of \cite{coc2} tells that sets satisfying
$\ufin(\rmO,\Gamma)$ and not containing perfect sets are perfectly meager.
In particular, $G$ is meager, and thus so is $L$, a contradiction.

(2) By Theorem 23 and Lemma 24 of \cite{o-bdd}, there is a Sierpi\'nski set\footnote{$S\sbst\R$ is \emph{Sierpi\'nski set}
if it is uncountable, but has countable intersection with each Lebesgue measure zero set.}
$S\sbst\R$ satisfying $\sone(\BG,\BG)^\uparrow$. (This result can alternatively and more easily be proved like Theorem \ref{OmOm},
by using \emph{random reals} instead of Cohen reals.)
Let $H=\<S\>$. By Theorem \ref{gengp}, $H$ satisfies $\sone(\BG,\BG)^\uparrow$.
As $\sone(\BG,\BG)$ implies $\sfin(\B,\B)$, $\ufin(\B,\BG)^\uparrow$ implies $\sfin(\B,\B)^\uparrow$, which is $\sfin(\BO,\BO)^\uparrow$.

It remains to prove that $H$ does not satisfy $\sone(\rmO,\rmO)$.
Assume otherwise.
By Remark \ref{BorHered}, As $H$ satisfies $\sone(\BG,\BG)$, it is a $\sigma$-space \cite{CBC}.
By Theorem \ref{ssCtoH}, its subset $S$ also satisfies $\sone(\rmO,\rmO)$, and thus has Lebesgue measure zero,
a contradiction.
\epf

\ed
\begin{thebibliography}{99}
\bibitem{FrGps}
T. Banakh, D. Repov\v{s}, and L. Zdomskyy,
\emph{$o$-Boundedness of free topological groups},
Topology and its Applications \textbf{157} (2010), 466--481.

\bibitem{barju}
T.\ Bartoszy\'nski and H.\ Judah,
Set Theory: On the structure of the real line, A.\ K.\ Peters,
Massachusetts: 1995.

\Pa{ideals}{T. Bartoszy\'nski and B. Tsaban}{Hereditary topological diagonalizations and the Menger-Hurewicz Conjectures}{Proceedings of the American Mathematical Society}{134}{2006}{605}{615}

\bibitem{BlassHBK}
A. Blass,
\emph{Combinatorial cardinal characteristics of the continuum},
in: \textbf{Handbook of Set Theory} (M.\ Foreman, A.\ Kanamori, and M.\ Magidor, eds.),
Kluwer Academic Publishers, Dordrecht, to appear.
\texttt{http://www.math.lsa.umich.edu/\~{}ablass/hbk.pdf}

\bibitem{Brendle96}
J.\ Brendle,
\emph{Generic constructions of small sets of reals},
Topology and its Applications \textbf{71} (1996), 125--147.

\bibitem{BH03}
L.\ Bukovsk\'y and J.\ Hale\v{s},
\emph{On Hurewicz properties},
Topology and its Applications \textbf{132} (2003), 71--79.

\bibitem{BRR91}
L. Bukovsk\'y, I. Rec\l{}aw, and M. Repick\'y,
\emph{Spaces not distinguishing pointwise and quasinormal convergence of real functions},
Topology and its Applications \textbf{41} (1991), 25--41.

\bibitem{BRR01}
L. Bukovsk\'y, I. Rec\l{}aw, M. Repick\'y,
\emph{Spaces not distinguishing convergences of real-valued functions},
Topology and its Applications \textbf{112} (2001), 13--40.

\bibitem{CPAgamma}
K. Ciesielski, A. Mill\'an, and J. Pawlikowski,
\emph{Uncountable $\gamma$-sets under axiom $\mathrm{CPA}^\mathrm{game} \mathrm{cube}$},
Fundamenta Mathematicae \textbf{176} (2003),  143--155.

\bibitem{DHM}
M. D\v{z}amonja, M. Hru\'sak, and J. Moore,
\emph{Parametrized $\lozenge$ principles},
Transactions of the American Mathematical Society \textbf{356} (2004), 2281--2306.

\bibitem{FM88}
D. Fremlin and A. Miller,
\emph{On some properties of Hurewicz, Menger and Rothberger},
Fundamenta Mathematica \textbf{129} (1988), 17--33.

\bibitem{GM}
F. Galvin and A. Miller,
\emph{$\gamma$-sets and other singular sets of real numbers},
Topology and its Applications \textbf{17} (1984), 145--155.

\bibitem{GN}
J.\ Gerlits and Zs.\ Nagy,
\emph{Some properties of $C(X)$, I},
Topology and its Applications \textbf{14} (1982), 151--161.

\bibitem{FUfin}
G. Gruenhage and P. Szeptycki,
\emph{Fr\'echet-Urysohn for finite sets},
Topology and its Applications \textbf{151} (2005) 238--259.

\bibitem{Jordan}
F. Jordan,
\emph{There are no hereditary productive $\gamma$-spaces},
Topology and its Applications \textbf{155} (2008), 1786--1791.

\bibitem{JSW}
H. Judah, S. Shelah, and H. Woodin,
\emph{The Borel conjecture},
Annals of Pure and Applied Logic \textbf{50} (1990),  255--269.

\bibitem{coc2}
W. Just, A. Miller, M. Scheepers, and P. Szeptycki,
\emph{The combinatorics of open covers II},
Topology and its Applications \textbf{73} (1996), 241--266.

\bibitem{coc7}
Lj.\ Ko\v{c}inac and M.\ Scheepers,
\emph{Combinatorics of open covers (VII): Groupability},
Fundamenta Mathematicae \textbf{179} (2003), 131--155.

\bibitem{KocSurv}
Lj.\ Ko\v{c}inac,
\emph{Selected results on selection principles},
in: \textbf{Proceedings of the 3rd Seminar on Geometry and Topology} (Sh.\ Rezapour, ed.),
July 15-17, Tabriz, Iran, 2004, 71--104.

\bibitem{MilNonGamma}
A. Miller,
\emph{A Nonhereditary Borel-cover $\gamma$-set},
Real Analysis Exchange \textbf{29} (2003/4), 601--606.

\bibitem{MillerBC}
A. Miller,
\emph{The $\gamma$-Borel conjecture},
Archive for Mathematical Logic \textbf{44} (2005), 425--434.
Appendix in \arx{math.LO/0312308}

\bibitem{MilRelGa}
A. Miller,
\emph{Cardinal characteristic for relative $\gamma$-sets},
Topology and its Applications \textbf{156} (2009), 872--878.

\bibitem{Sakai88}
M.\ Sakai,
\emph{Property $C''$ and function spaces},
Proceedings of the American Mathematical Society \textbf{104} (1988), 917--919.

\bibitem{coc1}
M.\ Scheepers,
\emph{Combinatorics of open covers I: Ramsey theory},
Topology and its Applications \textbf{69} (1996), 31--62.

\bibitem{alpha_i}
M.\ Scheepers,
\emph{$C_p(X)$ and Arhangel'ski\u{\i}'s $\alpha_i$ spaces},
Topology and its Applications \textbf{89} (1998), 265--275.

\bibitem{wqn}
M.\ Scheepers,
\emph{Sequential convergence in ${\sf C}_p(X)$ and a covering property},
East-West Journal of Mathematics \textbf{1} (1999),
207--214.

\bibitem{LecceSurvey}
M.\ Scheepers,
\emph{Selection principles and covering properties in topology},
Note di Matematica \textbf{22} (2003), 3--41.

\Pa{CBC}{M. Scheepers and B. Tsaban}{The combinatorics of Borel covers}{Topology and its Applications}{121}{2002}{357}{382}

\Pa{TkaTGII}{M. Tkachenko}{Topological groups for topologists: Part II}{Bolet\'\i{}n de la Sociedad Matem\'atica Mexicana. Tercera Serie}
{6}{2000}{1}{41}

\bibitem{Tod95}
S. Todor\v{c}evi\'c,
\emph{Aronszajn orderings},
Duro Kurepa memorial volume,
Publications de l'Institut Math\'ematique (new series) \textbf{57} (1995), 29--46.

\Bc{ict}{B. Tsaban}{Some new directions in infinite-combinatorial topology}{Set Theory}{J. Bagaria and S. Todor\v{c}evic, eds.}{Trends in Mathematics, Birkh\"auser}{2006}{225}{255}

\Pa{o-bdd}{B. Tsaban}{$o$-bounded groups and other topological groups with strong combinatorial properties}{Proceedings of the American Mathematical Society}{134}{2006}{881}{891}

\Bc{AddQuad}{B. Tsaban}{Additivity numbers of covering properties}{Selection Principles and Covering Properties in Topology}
{L. Ko\v{c}inac, ed.}{Quaderni di Matematica 18, Seconda Universita di Napoli, Caserta}{2006}{245}{282}

\Pa{capinf}{B. Tsaban}{A new selection principle}{Topology Proceedings}{31}{2007}{319}{329}

\Bc{OPiT}{B. Tsaban}{Selection Principles and special sets of reals}{Open Problems in Topology II}{E. Pearl, ed.}{Elsevier B.V.}{2007}{91}{108}

\Pa{SFT}{B. Tsaban and L. Zdomskyy}{Combinatorial images of sets of reals and semifilter trichotomy}{Journal of Symbolic Logic}{73}{2008}{1278}{1288}

\bibitem{hH}
B. Tsaban and L. Zdomskyy,
\emph{Hereditarily Hurewicz spaces and Arhangel'ski\u{\i} sheaf amalgamations},
Journal of the European Mathematical Society, to appear.

\end{thebibliography}
